\documentclass[reqno,12pt,a4paper]{amsart}

\usepackage{mathrsfs}
\usepackage{amssymb}
\usepackage{amsfonts}
\usepackage{latexsym}
\usepackage{amsthm}
\usepackage{enumerate}

\usepackage{graphicx}
\def\lb{\label}

\newcommand{\er}[1]{\textrm{(\ref{#1})}}

\begin{document}


\renewcommand{\theequation}{\arabic{section}.\arabic{equation}}
\theoremstyle{plain}
\newtheorem{theorem}{\bf Theorem}[section]
\newtheorem{lemma}[theorem]{\bf Lemma}
\newtheorem{corollary}[theorem]{\bf Corollary}
\newtheorem{proposition}[theorem]{\bf Proposition}
\newtheorem{definition}[theorem]{\bf Definition}
\newtheorem{remark}[theorem]{\it Remark}

\def\a{\alpha}  \def\cA{{\mathcal A}}     \def\bA{{\bf A}}  \def\mA{{\mathscr A}}
\def\b{\beta}   \def\cB{{\mathcal B}}     \def\bB{{\bf B}}  \def\mB{{\mathscr B}}
\def\g{\gamma}  \def\cC{{\mathcal C}}     \def\bC{{\bf C}}  \def\mC{{\mathscr C}}
\def\G{\Gamma}  \def\cD{{\mathcal D}}     \def\bD{{\bf D}}  \def\mD{{\mathscr D}}
\def\d{\delta}  \def\cE{{\mathcal E}}     \def\bE{{\bf E}}  \def\mE{{\mathscr E}}
\def\D{\Delta}  \def\cF{{\mathcal F}}     \def\bF{{\bf F}}  \def\mF{{\mathscr F}}
\def\c{\chi}    \def\cG{{\mathcal G}}     \def\bG{{\bf G}}  \def\mG{{\mathscr G}}
\def\z{\zeta}   \def\cH{{\mathcal H}}     \def\bH{{\bf H}}  \def\mH{{\mathscr H}}
\def\e{\eta}    \def\cI{{\mathcal I}}     \def\bI{{\bf I}}  \def\mI{{\mathscr I}}
\def\p{\psi}    \def\cJ{{\mathcal J}}     \def\bJ{{\bf J}}  \def\mJ{{\mathscr J}}
\def\vT{\Theta} \def\cK{{\mathcal K}}     \def\bK{{\bf K}}  \def\mK{{\mathscr K}}
\def\k{\kappa}  \def\cL{{\mathcal L}}     \def\bL{{\bf L}}  \def\mL{{\mathscr L}}
\def\l{\lambda} \def\cM{{\mathcal M}}     \def\bM{{\bf M}}  \def\mM{{\mathscr M}}
\def\L{\Lambda} \def\cN{{\mathcal N}}     \def\bN{{\bf N}}  \def\mN{{\mathscr N}}
\def\m{\mu}     \def\cO{{\mathcal O}}     \def\bO{{\bf O}}  \def\mO{{\mathscr O}}
\def\n{\nu}     \def\cP{{\mathcal P}}     \def\bP{{\bf P}}  \def\mP{{\mathscr P}}
\def\r{\rho}    \def\cQ{{\mathcal Q}}     \def\bQ{{\bf Q}}  \def\mQ{{\mathscr Q}}
\def\s{\sigma}  \def\cR{{\mathcal R}}     \def\bR{{\bf R}}  \def\mR{{\mathscr R}}
 \def\cS{{\mathcal S}}
\def\bS{{\bf S}}  \def\mS{{\mathscr S}}
\def\t{\tau}    \def\cT{{\mathcal T}}     \def\bT{{\bf T}}  \def\mT{{\mathscr T}}
\def\f{\phi}    \def\cU{{\mathcal U}}     \def\bU{{\bf U}}  \def\mU{{\mathscr U}}
\def\F{\Phi}    \def\cV{{\mathcal V}}     \def\bV{{\bf V}}  \def\mV{{\mathscr V}}
\def\P{\Psi}    \def\cW{{\mathcal W}}     \def\bW{{\bf W}}  \def\mW{{\mathscr W}}
\def\o{\omega}  \def\cX{{\mathcal X}}     \def\bX{{\bf X}}  \def\mX{{\mathscr X}}
\def\x{\xi}     \def\cY{{\mathcal Y}}     \def\bY{{\bf Y}}  \def\mY{{\mathscr Y}}
\def\X{\Xi}     \def\cZ{{\mathcal Z}}     \def\bZ{{\bf Z}}  \def\mZ{{\mathscr Z}}
\def\O{\Omega}


\newcommand{\gA}{\mathfrak{A}}
\newcommand{\gB}{\mathfrak{B}}
\newcommand{\gC}{\mathfrak{C}}
\newcommand{\gD}{\mathfrak{D}}
\newcommand{\gE}{\mathfrak{E}}
\newcommand{\gF}{\mathfrak{F}}
\newcommand{\gG}{\mathfrak{G}}
\newcommand{\gH}{\mathfrak{H}}
\newcommand{\gI}{\mathfrak{I}}
\newcommand{\gJ}{\mathfrak{J}}
\newcommand{\gK}{\mathfrak{K}}
\newcommand{\gL}{\mathfrak{L}}
\newcommand{\gM}{\mathfrak{M}}
\newcommand{\gN}{\mathfrak{N}}
\newcommand{\gO}{\mathfrak{O}}
\newcommand{\gP}{\mathfrak{P}}
\newcommand{\gR}{\mathfrak{R}}
\newcommand{\gS}{\mathfrak{S}}
\newcommand{\gT}{\mathfrak{T}}
\newcommand{\gU}{\mathfrak{U}}
\newcommand{\gV}{\mathfrak{V}}
\newcommand{\gW}{\mathfrak{W}}
\newcommand{\gX}{\mathfrak{X}}
\newcommand{\gY}{\mathfrak{Y}}
\newcommand{\gZ}{\mathfrak{Z}}

\def\ve{\varepsilon}   \def\vt{\vartheta}    \def\vp{\varphi}    \def\vk{\varkappa}

\def\Z{{\mathbb Z}}    \def\R{{\mathbb R}}   \def\C{{\mathbb C}}
\def\T{{\mathbb T}}    \def\N{{\mathbb N}}   \def\dD{{\mathbb D}}
\def\B{{\mathbb B}}    \def\dS{{\mathbb S}}


\def\la{\leftarrow}              \def\ra{\rightarrow}      \def\Ra{\Rightarrow}
\def\ua{\uparrow}                \def\da{\downarrow}
\def\lra{\leftrightarrow}        \def\Lra{\Leftrightarrow}


\def\lt{\biggl}                  \def\rt{\biggr}
\def\ol{\overline}               \def\wt{\widetilde}
\def\no{\noindent}


\let\ge\geqslant                 \let\le\leqslant
\def\lan{\langle}                \def\ran{\rangle}
\def\/{\over}                    \def\iy{\infty}
\def\sm{\setminus}               \def\es{\emptyset}
\def\ss{\subset}                 \def\ts{\times}
\def\pa{\partial}                \def\os{\oplus}
\def\om{\ominus}                 \def\ev{\equiv}
\def\iint{\int\!\!\!\int}        \def\iintt{\mathop{\int\!\!\int\!\!\dots\!\!\int}\limits}
\def\el2{\ell^{\,2}}             \def\1{1\!\!1}
\def\sh{\sharp}
\def\wh{\widehat}

\def\where{\mathop{\mathrm{where}}\nolimits}
\def\as{\mathop{\mathrm{as}}\nolimits}
\def\Area{\mathop{\mathrm{Area}}\nolimits}
\def\arg{\mathop{\mathrm{arg}}\nolimits}
\def\const{\mathop{\mathrm{const}}\nolimits}
\def\det{\mathop{\mathrm{det}}\nolimits}
\def\diag{\mathop{\mathrm{diag}}\nolimits}
\def\diam{\mathop{\mathrm{diam}}\nolimits}
\def\dim{\mathop{\mathrm{dim}}\nolimits}
\def\dist{\mathop{\mathrm{dist}}\nolimits}
\def\Im{\mathop{\mathrm{Im}}\nolimits}
\def\Iso{\mathop{\mathrm{Iso}}\nolimits}
\def\Ker{\mathop{\mathrm{Ker}}\nolimits}
\def\Lip{\mathop{\mathrm{Lip}}\nolimits}
\def\rank{\mathop{\mathrm{rank}}\limits}
\def\Ran{\mathop{\mathrm{Ran}}\nolimits}
\def\Re{\mathop{\mathrm{Re}}\nolimits}
\def\Res{\mathop{\mathrm{Res}}\nolimits}
\def\res{\mathop{\mathrm{res}}\limits}
\def\sign{\mathop{\mathrm{sign}}\nolimits}
\def\span{\mathop{\mathrm{span}}\nolimits}
\def\supp{\mathop{\mathrm{supp}}\nolimits}
\def\Tr{\mathop{\mathrm{Tr}}\nolimits}
\def\BBox{\hspace{1mm}\vrule height6pt width5.5pt depth0pt \hspace{6pt}}


\newcommand\nh[2]{\widehat{#1}\vphantom{#1}^{(#2)}}
\def\dia{\diamond}

\def\Oplus{\bigoplus\nolimits}



\def\qqq{\qquad}
\def\qq{\quad}
\let\ge\geqslant
\let\le\leqslant
\let\geq\geqslant
\let\leq\leqslant

\newcommand{\bea}{\begin{aligned}}
\newcommand{\ena}{\end{aligned}}

\newcommand{\ca}{\begin{cases}}
\newcommand{\ac}{\end{cases}}
\newcommand{\ma}{\begin{pmatrix}}
\newcommand{\am}{\end{pmatrix}}
\renewcommand{\[}{\begin{equation}}
\renewcommand{\]}{\end{equation}}
\def\bu{\bullet}

\title[{Trace formulas for Schr\"odinger operators on lattices}]
{Trace formulae for Schr\"odinger operators with complex-valued
potentials on cubic lattices}

\date{\today}

\author[Evgeny Korotyaev]{Evgeny Korotyaev}
\address{Saint-Petersburg State University, Universitetskaya nab. 7/9, St. Petersburg, 199034, Russia,
\ korotyaev@gmail.com, \ e.korotyaev@spbu.ru }

\author[Ari Laptev]{Ari Laptev}
\address{ Imperial College London, United Kingdom, \  a.laptev@imperial.ac.uk
}


\begin{abstract}
We consider  Schr\"odinger operators with complex decaying
potentials on the lattice. Using some classical results from Complex Analysis we  obtain some trace formulae and
using them estimate globally all zeros of the  Fredholm determinant in terms
of the potential.

\end{abstract}

\maketitle

\section{Introduction}

\noindent
Let us consider the Schr\"odinger operator $H$ acting in
$\ell^2(\Z^{d}),d\ge 3$ and given by
$$
{H}=H_0+V, \qqq H_0=\D,
$$
where ${\D}$ is the discrete Laplacian on $\Z^d$ given by
$$
\big(\D f \big)(n)=\frac{1}{2}\sum_{j=1}^{d}\big(f(n+ e_{j}) + f(n-
e_{j})\big),\qqq n=(n_j)_{j=1}^d\in \Z^d,
$$
for $f =(f_n)_{n\in\Z^d} \in \ell^{2}(\Z^d)$. Here
$
e_{1} = (1,0,\cdots,0), \cdots, e_{d} =
(0,\cdots,0,1)
$
 is the standard basis of $\Z^d$.
The operator  $V = (V_n)_{n\in\Z^d}$, $V_n\in\C$, is a complex
potential given by
\begin{equation*}
(Vf)(n)=V_nf_n, \qqq n\in \Z^d.
\end{equation*}
We assume that the potential $V$ satisfies  the following condition:
\begin{equation}
\lb{Vc}
V\in \ell^{2/3}(\Bbb Z^d).
\end{equation}
Note that the condition \eqref{Vc} implies that $V$ can be factorised as
\begin{equation} \lb{Vfact}
V = V_1V_2, \qquad {\rm where}\quad V_1\in \ell^1(\Bbb Z^d), \, V_2 \in \ell ^2(\Bbb Z^d),
\end{equation}
with $V_1=|V |^{2/3-1} V$ and $V_2=|V|^{1/3}$.

\noindent
Here  $\ell^q(\Z^{d}), q>0$ is the space of sequences
 $f=(f_n)_{n\in \Z^d}$ such that  $\|f\|_{q}<\iy$, where
$$
\begin{aligned}
\|f\|_{q}=\|f\|_{\ell^q(\Z^{d})} =
\begin{cases}   \sup_{n\in \Z^d}|f_n|,\qq &
\ q=\iy, \\
\big(\sum_{n\in \Z^d}|f_n|^q\big)^{1\/q},\qq & \ q\in (0,\iy).
\end{cases}
\end{aligned}
$$
Note that $\ell^q(\Z^{d}), q\ge 1$ is the Banach space equipped with the norm $\|\cdot\|_{q}$.
It is well-known that the
spectrum of the Laplacian $\D$ is absolutely continuous and equals
$$
\s(\D)=\s_{\textup{ac}}(\D)=[-d,d].
$$
It is also well known that if $V$ satisfies \eqref{Vc}, the essential spectrum of the Schr\"odinger
operator $H$ on $\ell^2(\Z^d)$ is
$$
\s_{\textup{ess}}(H)=[-d,d].
$$
However, this condition does not exclude appearance of the singular continuous spectrum on the interval $[-d,d]$.
Our main goal is to find new trace formulae for the operator $H$ with
complex potentials $V$ and to use these  formulae for some estimates of complex eigenvalues in terms of potentials.

Note that some of the results obtained in this paper are new even in the case of real-valued potentials due to the presence of the measure $\sigma$ (see Theorem \ref{T3}) appearing in the canonical factorisation of the respective Fredholm determinants. Non-triviality of such a measure is due to the weak condition \eqref{Vc} on the potential $V$. We believe that it would be interesting to study the relation between properties of $V$v and $\sigma$.

\medskip
Recently uniform bounds on eigenvalues of Schr\"odinger operators in $\Bbb R^d$ with complex-valued potentials decaying at infinity attracted attention of many specialists. We refer to \cite{Da} for a review of the state of the art of non-selfadjoint Schr\"odinger operators and for motivations and applications. Bounds on single eigenvalues were proved, for instance, in \cite{AAD,DN,FrLaSe,Fr} and bounds on sums of powers of eigenvalues were found in \cite{FrLaLiSe,LaSa,DeHaKa0,DeHaKa,BGK,FrSa,Fr3}. The latter bounds generalise the Lieb--Thirring bounds \cite{LiTh} to the non-selfadjoint setting.
Note that in \cite{FrSa} (Theorem 16) the authors obtained estimates on the sum of the distances between the complex eigenvalues and the continuous spectrum $[0,\infty)$ in terms of $L^p$-norms of the potentials.
Note that almost no results are known on the number of eigenvalues of Schr\"odinger operators with complex potentials. We referee here to a recent paper \cite{FLS} where the authors discussed this problem in details in odd dimensions.

\medskip
For the discrete Schr\"odinger operators most of the results were obtained in the self-adjoint case, see, for example,  \cite{T89} (for the $\Z^1$ case).  Schr\"odinger operators with decreasing potentials
on the lattice $\Z^d$ have been considered by Boutet de Monvel-Sahbani
\cite{BS99}, Isozaki-Korotyaev \cite{IK12}, Kopylova \cite{Kop10},
Rosenblum-Solomjak \cite{RoS09}, Shaban-Vainberg \cite{SV01} and see
references therein. Ando \cite{A12} studied the inverse spectral
theory for the discrete Schr\"odinger operators with finitely
supported potentials on the hexagonal lattice.
Scattering on periodic metric graphs  $\Bbb Z^d$ was considered by Korotyaev-Saburova \cite{KS}.

\noindent
Isozaki and Morioka (see Theorem 2.1. in \cite{IM14}) proved that if
the potential $V$ is real and compactly supported, then the
point-spectrum of $H$ on the interval $(-d,d)$ is absent. Note that in [10] the author gave an example of embedded eigenvalue at the
endpoint $\{\pm d\}$.

\medskip
\noindent
In this paper we use classical results from Complex Analysis that lead us to a new class of trace formula for the spectrum of discrete multi-dimensional Scr\"odinger operators with complex-valued potentials. In particular, we consider a so-called canonical factorisation of analytic functions from Hardy spaced via its inner and outer factors, see Section 6.
Such factorisations allied for Fredholm determinants allow us to obtain trace formula that lead to some inequalities on the complex spectrum in terms of the $L^{2/3}$ norm of the potential function.
Note also that in the case $d=3$ we use a delicate uniform inequality for Bessel's functions obtained in Lemma \ref{case3}.


\section{Some notations and statements of main results}

We denote by  $\dD_r(z_0)\subset\C$   the disc with radius $r>0$ and
center $z_0\in \C$
$$
\dD_r(z_0)=\{z\in \C:|z-z_0|<r\},
$$
and abbreviate $\dD_r=\dD_r(0)$ and $\dD=\dD_1$. Let also $\Bbb T = \partial \Bbb D$. It is convenient to introduce a new spectral
variable $z\in \dD$ by
\begin{equation}\label{lambda}
\l=\l(z)={d\/2}\rt(z+{1\/z}\rt)\in \L=\C\sm [-d,d] ,\qqq z\in \dD.
\end{equation}
 The function $\l(z)$ has  the following properties:

\medskip
{\it $\bu$   $z\to \l(z)$ is a conformal mapping from $\dD$
onto the spectral domain $\L$.

$\bu$  $\l(\dD)=\L=\C\sm [-d,d] $  and   $\l(\dD\cap \C_\mp)=\C_\pm
$.

$\bu$  $\L$ is the cut domain with the cut $[-d,d]$, having the
upper side $[-d,d]+i0$ and the lower side $[-d,d]-i0$.
$\l(z)$ maps the boundary: the upper semi-circle onto the lower side
$[-d,d]-i0$ and the lower semi-circle onto the upper side
$[-d,d]+i0$.

$\bu$  $\l(z)$ maps  $z=0$ to
$\l=\iy$.

$\bu$   The inverse mapping $z(\cdot ): \L\to \dD$ is given by
$$
\begin{aligned}
z={1\/d}\rt(\l-\sqrt{\l^2-d^2}\rt),\qqq \l\in \L,\\
z={d\/2\l}+{O(1)\/\l^3},\qqq as \qq |\l|\to \iy.
\end{aligned}
$$
}

Next we introduce the Hardy space $ \mH_p=\mH_p(\dD)$. Let $F$ be analytic in
$\dD$. For $0<p\le \iy$  we say $F$ belongs the Hardy space $ \mH_p$
if $F$ satisfies $\|F\|_{\mH_p}<\iy$, where $\|F\|_{\mH_p}$ is given
by
$$
\|F\|_{\mH_p}=
\begin{cases}
\sup_{r\in (0,1)}\rt({1\/2\pi}\int_\T
|F(re^{i\vt})|^pd\vt\rt)^{1\/p}, &
if  \qqq 0< p<\iy,\\
 \sup_{z\in \dD}|F(z)|, & if \qqq p=\iy.
 \end{cases}
$$
Let  $\cB$ denote the class of bounded operators and  $\cB_1$ and
$\cB_2$ be the trace and the Hilbert-Schmidt class equipped with the
norm $\|\cdot \|_{\cB_1}$ and $ \|\cdot \|_{\cB_2}$ respectively.

\medskip
Denote by $D(z), z\in \dD$  the determinant
$$
D(z)=\det \left(I+VR_0(\l(z)\right), \qqq z\in \dD,
$$
where $R_0(\l)=(H_0-\l)^{-1}, \l\in \L$.
The determinant  $D(z), z\in \dD$, is well
defined for  $V\in \cB_1$ and  if $\l_0\in \L$
is an eigenvalue of $H$, then $z_0=z(\l_0)\in \dD$ is a zero of $D$
with the same multiplicity.

\begin{theorem}
\lb{T1} Let a potential $V$ satisfy \er{Vc}. Then the determinant
$D(z)=\det (I+VR_0(\l(z)), \ z\in \dD$, is analytic in $\Bbb D$. It has $N\le \iy$ zeros $\{z_j\}_{j=1}^N$, such that

\[
\lb{zD}
\begin{aligned}
& 0<r_0=|z_1|\le |z_2|\le ...\le |z_j|\le |z_{j+1}\le
|z_{j+2}|\le .... ,\\
& where \qqq r_0=\inf |z_j|>0.
\end{aligned}
\]
Moreover, it satisfies
\begin{equation}
\label{D1}
\begin{aligned}
\|D\|_{\mH_\iy}\le e^{C\|V\|_{2/3}},
\end{aligned}
\end{equation}
where the constant $C$ depends only on $d$.

\medskip
\noindent
Furthermore, the function $\log D(z)$ whose branch is defined by $\log D(0)=0$,  is
analytic in the disk $\dD_{r_0}$ with the radius  $r_0>0$ defined by
\er{zD} and it has the Taylor series as $|z|<r_0$:
\[
\label{D3} \log D(z)=-c_1z-c_2z^2-c_3z^3-c_4z^4 -.....
\]
 where
\[
\label{D4} c_1=d_1a,\qqq c_2=d_2a^2, \qqq c_3=d_3a^3-c_1, \qqq
c_4=d_4a^4-c_2,....
\]
\begin{equation}
\label{cD4}
\begin{aligned}
d_1=\Tr V,\qqq d_2= \Tr\,V^2,\qqq
d_3=\Tr\,\big(V^3+(3d/2)V\big),....,
\end{aligned}
\end{equation}
and where $a={2\/d}$.
\end{theorem}

\medskip
\noindent
Define the Blaschke product $B(z), z\in \dD$ by
\[
\lb{B2}
\begin{aligned}
& B(z)=\prod_{j=1}^N {|z_j|\/z_j}{(z_j-z)\/(1-\ol z_j z)},\qqq
 & if \qqq N\ge 1,\\
&B=1, \qqq & if \qqq  N=0.
\end{aligned}
\]


\begin{theorem}\lb{T2}
Let a potential $V$ satisfy \er{Vc} and let $N\ge 1$. Then the zeros
$\{z_j\}$ of $D$ in the disk $\dD$ labeled by \er{zD}  satisfy
\[
\lb{B1} \sum _{j=1}^N (1-|z_j|)<\iy.
\]
Moreover, the Blaschke product $B(z), z\in \dD$ given by \er{B2}
converges absolutely for $\{|z|<1\}$ and  satisfies

i)  $B\in \mH_\iy$ with $\|B\|_{\mH_\iy}\le 1$,
\[
\lb{B3} \lim_{r\to 1}|B(re^{i\vt})|=|B(e^{i\vt})|=1 \qq   {\rm for \ almost \ all}Ê\qq \vt\in \T,
\]
and
\[
\lb{B4} \lim _{r\to 1}\int_0^{2\pi}\log |B(re^{i\vt})|d\vt=0.
\]

ii) The determinant $D$ has the factorization in the disc $\dD$:
\begin{equation*}
D=BD_B,
\end{equation*}
 where $D_B$ is analytic in the unit disc $\dD$
and has not zeros in $\dD$.

 iii) The Blaschke product $B$ has the Taylor series at $z=0$:
\[
\begin{aligned}
\lb{B6}
\log  B(z)=B_0-B_1z-B_2z^2-... \qqq as \qqq z\to 0,
\end{aligned}
\]
where $B_n$ satisfy
\begin{equation*}
\begin{aligned}
& B_0=\log  B(0)<0,\qqq B_1=\sum_{j=1}^N\rt({1\/z_j}-\ol z_j \rt),...,
\qqq B_n={1\/n}\sum_{j=1}^N\rt({1\/z_j^n}-\ol z_j^n \rt),....\\
& |B_n|\le {2\/r_0^n}\sum _{j=1}^N (1-|z_j|).
\end{aligned}
\end{equation*}
\end{theorem}


\bigskip
\noindent
The next statement describes the canonical representation of the determinant $D(z)$.

\begin{theorem}
\lb{T3} Let a potential $V$ satisfy \er{Vc}. Then

i)   There exists a singular measure $\s\ge 0$ on $[-\pi,\pi]$, such
that the determinant $D$ has a canonical factorization for all
$|z|<1$ given by
\[
\lb{cfD}
\begin{aligned}
& D(z)=B(z)e^{-K_\s (z)}e^{K_D(z)},\\
& K_\s(z)={1\/2\pi}\int_{-\pi}^{\pi}{e^{it}+z\/e^{it}-z}d\s(t),\\
& K_D(z)= {1\/2\pi}\int_{-\pi}^{\pi}{e^{it}+z\/e^{it}-z}\log
|D(e^{it})|dt,
 \end{aligned}
\]
 where $\log |D(e^{it}) |\in L^1(-\pi,\pi)$.

ii) The measure $\s$ satisfies
\[
\supp \s\ss \{t\in [-\pi,\pi]: D(e^{it})=0\}.
\]
\end{theorem}

\noindent
{\bf Remarks.}

\smallskip
\noindent
1) For the canonical factorisation of analytic functions see, for example, \cite{Koo98}.

\smallskip
\noindent
2)
Note that for $D_{in}(z)$ defined by $D_{in}(z)= B(z) e^{-K_\sigma(z)}$,
we have  $| D_{in}(z)|\le 1$, since $d\s\ge
0$ and  $\Re {e^{it}+z\/e^{it}-z}\ge 0$ for all $(t,z)\in
\T\ts\dD$.

\smallskip
\noindent
3) The closure of the set $\{z_j\}\cup \supp \s$ is called the
spectrum of the inner  function $ D_{in}$.

\smallskip
\noindent
4) $D_B={D\/B}$  has no zeros in the disk $\dD$ and
satisfies
$$
\log D_B(z)={1\/2\pi}\int_{-\pi}^{\pi}{e^{it}+z\/e^{it}-z}d\m(t),
$$
where the measure $\m$ equals
$$
 d\m(t)=\log
|D(e^{it})|dt-d\s(t).
$$

\bigskip
\begin{theorem}
\lb{T4} {\bf (Trace formulae.)}  Let $V$ satisfy
\er{Vc}. Then the following identities hold
\[
\lb{tr0} {\s(\T)\/2\pi}-B_0={1\/2\pi}\int_{-\pi}^{\pi}\log |D(e^{it})|dt\ge 0,
\]
\[
\lb{tr1} -D_n+B_n={1\/\pi}\int_\T e^{-int}d\m(t),\qqq n=1,2,....
\]
where $B_0=\log B(0)=\log \left(\prod_{j=1}^N |z_j|\right)<0$ and $B_n$ are given by
\er{B6}. In particular,
\[
\lb{tr4} \sum_{j=1}^N\rt({1\/z_j}-\ol z_j \rt)={2\/d}\Tr \,
V+{1\/\pi}\int_\T e^{-it}d\m(t),
\]
\[
\lb{tr55} \sum_{j=1}^N\rt({1\/z_j^2}-\ol z_j^2 \rt)={4\/d^2}\Tr \,
V^2+{1\/\pi}\int_\T e^{-i2t}d\m(t),
\]
and
\[
\lb{t52}
\begin{aligned}
 \sum_{j=1}^N \Im\l_j=\Tr \Im V-{d\/2\pi}\int_\T \sin t\,
d\m(t),
\\
\sum_{j=1}^N \Re\sqrt{\l_j^2-d^2}=\Tr \Re V+{d\/2\pi}\int_\T \cos
t\,d\m(t).
\end{aligned}
\]

\end{theorem}


\bigskip
\begin{theorem}
\lb{T5} Let  $V$ satisfy \er{Vc}. Then we have the following
estimates:
\[
\lb{t51} \sum (1-|z_j|)\le -B_0\le C(d)\|V\|_{2/3}-{\s(\T)\/2\pi}.
\]
and if $\Im V\ge 0$, then
\[
\lb{t51x}
 \sum_{j=1}^N \Im\l_j\le \Tr \Im V+C(d)\|V\|_{2/3},
\]
and if $ V\ge 0$, then
\[
\lb{t51xx}
 \sum_{j=1}^N \sqrt{\l_j^2-d^2}\le \Tr V+C(d)\|V\|_{2/3},
\]

\end{theorem}


\noindent
{\bf Remark.}

\noindent
Note that some of the results stated in Theorems \ref{T4} and \ref{T5} are new even for real-valued  potentials, see Section 5.


\section {Determinants  }
\setcounter{equation}{0}

\subsection {Properties of the Laplacian } One may diagonalize the
discrete Laplacian, using the (unitary) Fourier transform
 $\F\colon \ell^2(\Z^d)\to L^2(\Bbb S^d)$, where $\Bbb S=\R/(2\pi \Z)$. It is defined by
 $$
 (\F f)(k)=\wh f(k)={1\/(2\pi)^{{d\/2}}}\sum_{n\in
 \Z^d} f_ne^{i(n,k)},\qq \textup{where} \qq
 k=(k_j)_{j=1}^d\in \Bbb S^d.
 $$
 Here $(\cdot,\cdot)$ is the scalar product in $\R^d$.
In the so-called momentum representation of the operator $H$, we have:
$$
 \F H \F^*=\wh \D +\wh V.
$$
The Laplacian is transformed into the multiplication operator
$$
(\wh \D \wh f)(k)=h(k)\wh   f(k),\qqq h(k)=\sum_1^d \cos k_j,\qqq
k\in \Bbb S^d,
$$
and the potential $V$ becomes a convolution operator
$$
(\wh V\wh f)(k)={1\/(2\pi)^{d\/2}}\int_{\Bbb S^d} \wh V(k-k')\wh f(k')dk',
$$
where
$$
\wh V(k)={1\/(2\pi)^{d\/2}}\sum_{n\in \Z^d}V_ne^{i(n,k)},\qqq
V_n={1\/(2\pi)^{d\/2}}\int_{\Bbb S^d} \wh V(k)e^{-i(n,k)}dk.
$$

\subsection{Trace class operators}
Here for the sake of completeness we give some standard facts from Operator Theory in Hibert Spaces..

\noindent
Let $\mathcal H$ be a Hilbert space endowed with inner product $(\,
, \, )$ and norm $\|\cdot\|$. Let $\cB_1$ be the set of all trace
class operators on $\mathcal H$ equipped with the trace norm
$\|\cdot\|_{\cB_1}$. Let us recall some well-known facts.

$\bu$ Let $A, B\in \cB$ and $AB, BA\in \cB_1$. Then
\begin{equation*}
\label{AB} {\rm Tr}\, AB={\rm Tr}\, BA,
\end{equation*}
\begin{equation*}
\label{1+AB} \det (I+ AB)=\det (I+BA).
\end{equation*}
\begin{equation*}
\label{DA1} |\det (I+ A)|\le e^{\|A\|_{\cB_1}}.
\end{equation*}
\begin{equation*}
\label{DA1x} |\det (I+ A)-\det (I+ B)|\le \|A-B\|_{\cB_1}
e^{1+\|A\|_{\cB_1}+\|B\|_{\cB_1}}.
\end{equation*}
Moreover,  $I+ A$ is invertible if and only if $\det (I+ A)\ne 0$.

%
%
%
%

$\bu$  Suppose for a domain $\mD \subset {\C}$, the function
$\O(\cdot)-I: \O\to \cB_1 $ is analytic and  invertible
for any $z\in D$. Then for  $F(z)=\det \Omega (z)$ we have
\begin{equation*}
 F'(z)= F(z){\rm Tr}\,\left(\O(z)^{-1}\O'(z)\right).
\end{equation*}

$\bu$  Recall that for $K \in \cB_1$ and $z \in {\C}$, the following
identity holds true:
\begin{equation*}
\det\,(I - zK) = \exp\left(- \int_0^z{\rm Tr}\,
\big(K(1 - sK)^{-1}\big)ds\right) \label{S6Detdefine}
\end{equation*}
(see e.g. \cite{GK}, p.167, or \cite{RS78}, p.331).

\bigskip


\subsection{Fredholm determinant}
We recall here results from \cite{IK12} about the asymptotics of the
determinant $\cD(\l)=\det (I+VR_0(\lambda))$ as $|\l|\to \iy$.

\begin{lemma}
\label{TaD1} Let  $V\in \ell^1(\Z^d)$. Then the determinant
$\cD(\l)=\det (I+VR_0(\lambda))$ is analytic in  $\L=\C\sm [-d,d]$
and satisfies
\begin{equation*}
\cD(\l)=1+O(1/\l) \quad as \quad |\l|\to {\infty},
\end{equation*}
uniformly in $\arg \l\in [0,2\pi]$, and
\begin{equation*}
 \log \cD(\l) = - \sum_{n=1}^{\infty}\frac{(-1)^n}{n}{\rm
Tr}\,\left(VR_0(\l)\right)^n,
\end{equation*}
\begin{equation}
\label{aD3} \log \cD(\l) =-\sum _{n \geq 1}\frac{d_n}{n\l^n},\quad
 d_n={\rm Tr}\,(H^n-H_0^n),
\end{equation}
where the right-hand side is absolutely convergent for $|\lambda| > r_1$, $r_1 >0$ being a sufficiently large constant.  In particular,
\begin{equation}
\label{aD4}
\begin{aligned}
&d_1={\Tr} \, V,\\
&d_2={\Tr}\,V^2,\\
& d_3={\Tr}\,\big(V^3+6d\tau^2V\big),\\
& d_4={\rm Tr}\, \big(V^4+8d\tau^2V^2+2\tau^2(V_\D)V\big), \dots,
\end{aligned}
\end{equation}
where $V_\D=\sum_{i=1}^d(S_jVS_j^*+S_j^*VS_j)$ and
$(S_jf)(n)=f(n+e_j)$ and $\t={1\/2}$.
\end{lemma}

\bigskip
\noindent
Recall the conformal mapping $\l(\cdot): \dD\to \L$ is
given by $ \l(z)={d\/2}\big(z+{1\/z} \big),\ |z|<1$, and note that
$|\l|\to \iy$ iff $z\to 0$. We consider the operator-valued function
$Y(\l(z)), \l\in \L$, defined by
$$
Y(\l(z)) = V_2 X (\l(z)), \qquad {\rm where}\qquad X (\l(z)) = |V_1|^{1/2} R_0 (\l(z)) |V_1|^{1/2} V_1|V_1|^{-1},
$$
and where $V_1$ and $V_2$ are defined in \er{Vfact}.


\begin{theorem}
\lb{T2x}
Let $V$ satisfy \er{Vc}. Then the
operator-valued function $Y(\l(z)): \dD\to \cB_1$ is analytic in
the unit disc $\dD$ and satisfies
\begin{equation*}
\begin{aligned}
\|Y(\l(z))\|_{\cB_1}\le C(d) \|V\|_{2/3},\qqq \forall \
z\in \dD.
\end{aligned}
\end{equation*}
Moreover, the function $D(z), z\in \dD$ belongs
to $\mH_\iy$ and
\[
\label{det1}
\|D(\cdot)\|_{\mH_\iy}\le e^{C(d) \|V\|_{2/3}}.
\]

\end{theorem}

\medskip
\noindent
{\bf Proof.}
The operator $V_2$ belongs to $\cB_2$ and due to Theorem \ref{TApp} (see Appendix 2) the operator-function $X(\cdot):\L\to \cB_2$
satisfies the inequality
\begin{equation*}
\begin{aligned}
\|X(\l)\|_{\cB_2}\le C(d) \|V_1\|_2,\qqq \forall\
\l \in \Lambda =  \C\sm [-d,d],
\end{aligned}
\end{equation*}
Thus  the operator-valued function $Y(\l(z))=V_2\, X(\l(z)):
\dD\to \cB_1$ is of trace class. Moreover, the function
$D(z), z\in \dD$, belongs to $\mH_\iy$ and due to \er{DA1} it
satisfies \er{det1}. \BBox

\bigskip
The function $D(z)=\cD(\l(z))$ is analytic in $\dD$ with the zeros
given by
$$
z_j=z(\l_j), \qqq j=1,2,..., N,
$$
where $\l_j$ are zeros (counting with multiplicity) of $\cD(\l)$ in
$\L=\C\sm [-d,d]$, i.e., eigenvalues of $H$.


\begin{lemma}
\label{TaD2} Let a potential $V\in \ell^1(\Z^d)$. Then
$\log D(z)$ is analytic in $\dD_{r_0}$ defined by $\log D(0)=0$,
where $r_0$ is given by \er{zD}, and has the following Taylor series
\[
\label{aaD1} \log D(z)=-c_1z-c_2z^2-c_3z^3-c_4z^4-......, \qqq \as
\qq |z|<r_0,
\]
 and
\[
\label{aaD2} c_1=d_1a,\qqq c_2=d_2a^2, \qqq c_3=d_3a^3-c_1, \qqq
c_4=d_4a^4-c_2,....
\]
where $a={2\/d}$ and the coefficients $d_j$ are given by \er{aD4}.
\end{lemma}


\bigskip
\noindent
{\bf Proof.}   We have
$$
{1\/\l}={az\/1+z^2}=a(z-z^3+O(z^5)),\qqq {1\/\l^2}=a^2(z^2-z^4+O(z^6)),
\qqq {1\/\l^3}=a^3z^3+O(z^5)
$$
as $z\to 0$. Substituting these asymptotics into \er{aD3} we obtain
\er{aaD1} and \er{aaD2}. \BBox


\bigskip
\section {Proof of the main results}
\setcounter{equation}{0}

\noindent
We are ready to prove main results.

\bigskip
\noindent
{\bf Proof of Theorem \ref{T1}.}  Let  $V$ satisfy
\er{Vc}. Then by Theoem \ref{T2x}, the determinant $D(z), z\in
\dD$, is analytic and $D\in \mH_\iy$.
Moreover, Lemma \ref{TaD1} gives  \er{D1} and Lemma \ref{TaD2} gives
\er{D3}-\er{cD4}. \BBox

\bigskip
\noindent
{\bf Proof of Theorem \ref{T2}.}
Due to Theorem \ref{T1} the determinant $D(z)$ is analytic in  $\dD$.
Then  Theorem \ref{TA1} (see Appendix 1) yields
$$
\cZ_D:=\sum_{j=1}^\iy(1-|z_j|)<\iy
$$
 and
the Blaschke product $B(z)$ given by
$$
 B(z)=\prod_{j=1}^N {|z_j|\/z_j}{z_j-z\/1-\ol z_j z}, \qqq
 z\in \dD,
$$
 converges absolutely for $\{|z|<1\}$.  We have
$D(z)=B(z)D_B(z)$, where $D_B$ is analytic in the unit disc $\dD$
and has no zeros in $\dD$. Thus we have proved  ii).

i) Lemma \ref{TF2} gives \er{B3} and \er{B4}.

iii)  For small a sufficiently small $z$ and for $t=z_j\in \dD$ for some
$j$ we have the following identity:
$$
\log {|t|\/t}{t-z\/1-\ol t z}=\log |t|+ \log \rt(1-{z\/t}\rt)-\log
(1-\ol t z)
=\log |t|-\sum_{n\ge1}\rt({1\/t^n}-\ol t^n\rt){z^n\/n}.
$$
Besides,
$$
\begin{aligned}
&  |1-|t|^n|\le n|1-|t||,\\
 &\big|t^{-n}-\ol t^n \big|\le \big|1-t^n| +\big|1-t^{-n} \big|\le
|1-t^n|\rt(1+{1\/|t|^n} \rt)\le
 |1-|t|^n|{2\/r_0^n}\le  |1-|t||{2n\/r_0^n},
\end{aligned}
$$
where $r_0=\inf |z_j|>0$.
 This yields
\[
\begin{aligned}
& \log  B(z)=\sum_{j=1}^N\log  {|z_j|\/z_j}{z_j-z\/1-\ol z_j z}
=\sum_{j=1}^N\rt( \log |z_j|+ \log \big(1-(z/z_j)\big)-\log (1-\ol z_j
z)\rt)\\
&=\sum_{j=1}^N\log |z_j|-\sum_{n=1}^\iy\sum_{j=1}^N\rt({1\/z_j^n}-\ol
z_j^n \rt){z^n\/n}=\log B(0)-b(z),\\
& b(z)=\sum_{n=1}^N\sum_{j=1}^\iy\rt({1\/z_j^n}-\ol z_j^n
\rt){z^n\/n}=\sum_{n=1}^N z^nB_n,\qqq
B_n={1\/n}\sum_{j=1}^N\rt({1\/z_j^n}-\ol z_j^n \rt),
\end{aligned}
\]
where  the function $b$ is analytic in the disk $\{|z|<{r_0\/2}\}$  and $B_n$ satisfy
$$
\begin{aligned}
 |B_n|\le {1\/n}\sum_{j=1}^N\rt|{1\/z_j^n}-\ol z_j^n \rt| \le
{2\/r_0^n}\sum_{j=1}^N  |1-|z_j||={2\/r_0^n}\cZ_D,
 \end{aligned}
$$
where $\cZ_D=\sum_{j=1}^\iy(1-|z_j|)$.
 Thus
$$
|b(z)|\le \sum_{n=1}^\iy |B_n|{|z|^n}\le
2\cZ_D\sum_{n=1}^\iy{|z|^n\/r_0^n}={2\cZ_D\/1-{|z|\/r_0}}.
$$
\BBox

\bigskip
\noindent
{\bf Proof of Theorem \ref{T3}.}

 i) Theorem \ref{T1} implies $D\in \mH_\iy$.
Therefore the  canonical representation \er{cfD} follows from Lemma \ref{Tft}.

ii) The relation \er{meraze} gives the proof of ii). \BBox


\bigskip
\noindent
{\bf Proof of Theorem \ref{T4}.}  By using Lemma
\ref{TAt}, \er{D1}-\er{D4} and Proposition \ref{T2} we obtain
identities \er{tr0}-\er{tr55}.

We have the following identities for $z\in \dD$ and $\z={\l\/d}\in
\L_1$:
\[
\label{Et1}
\begin{aligned}
2\z=z+{1\/z},\qqq z=\z-\sqrt{\z^2-1},\qqq z-{1\/z}=-2\sqrt{\z^2-1}.
\end{aligned}
\]
These identities yield
\[
\label{Et2}
\begin{aligned}
 \ol z-{1\/z}=z+\ol z-2\z=2\Re z-2\z,\\
\ol z-{1\/z}=\ol z-z-2\sqrt{\z^2-1}=-2i\Im z-2\sqrt{\z^2-1}.
 \end{aligned}
\]
Then we get
$$
\begin{aligned}
{2\/d}\Tr \Im
V+\Im  \sum_{j=1}^N\rt(\ol z_j -{1\/z_j}\rt)={1\/\pi}\int_\T \sin t\, d\m(t),\\
\Im  \rt(\ol z_j -{1\/z_j}\rt)= -2\Im\z_j= -{2\/d}\Im\l_j,\\
\Re  \rt(\ol z_j -{1\/z_j}\rt)= 2\Re(z_j-\z_j)=-2\Re\sqrt{\z_j^2-1}
 \end{aligned}
$$
and thus
$$
\sum_{j=1}^N \Im\l_j=\Tr \Im V-{d\/2\pi}\int_\T  \sin t\,d \m(t),
$$
$$
\sum_{j=1}^N \Re\sqrt{\l_j^2-d^2}=\Tr \Re V+{d\/2\pi}\int_\T \cos
t\,d\m(t),
$$
\BBox

\bigskip
\noindent {\bf Proof of Theorem \ref{T5}.} The simple inequality
$1-x\le -\log x$ for $\forall \ x\in (0,1]$, implies
$-B_0=-B(0)=-\sum \log |z_j|\ge \sum (1- |z_j|)$. Then substituting
the last estimate and the estimate \er{D1} into the first trace
formula \er{tr0}  we obtain \er{t51}.

In order to determine the next two estimates we use the trace
formula \er{tr0}. Let $\Im V\ge 0$. Then $\Im \l_j\ge 0$ and the
estimates  \er{D1} and \er{t51} and the second trace formula
\er{t52} imply
$$
\begin{aligned}
 \sum_{j=1}^N \Im\l_j-\Tr \Im V=-{d\/2\pi}\int_\T \sin t\,
d\m(t)\\
\le {d\/2\pi}\int_\T (C\|V\|_{2\/3}dt+ d\s(t)) \le C(d)\|V\|_{2\/3},
\end{aligned}
$$
which yields \er{t51x}. Similar arguments give \er{t51xx}. \BBox


\section {Schr\"odinger operators with real potentials}
\setcounter{equation}{0}

Consider Schr\"odinger operators $H=-\D+V$, where the potential $V$
is real and satisfies condition \er{Vc} . The spectrum of $H$ has the
form
$$
\s(H)=\s_{ac}(H)\cup \s_{sc}(H)\cup \s_p(H)\cup \s_{dis}(H),\qq
\s_{ac}(H)=[-d,d],
$$
where
$$
  \s_p(H)\ss [-d,d], \qq \s_{dis}(H)\ss \R\sm[-d,d].
$$
Note that each eigenvalues of $H$ has a finite multiplicity.

\subsection{Discrete spectrum} The discrete eigenvalues of
the operator $H$ are real and belong to the set $\R\sm[-d,d]$. Let
they are labeled by
$$
\dots \le \lambda_{-2}   \le \lambda_{-1}<-d < d< \lambda_1\le \lambda_2 \le \dots
$$
The corresponding point from  $z_j\in\dD$ is real and satisfy
$$
\l_j={d\/2}\rt(z_j+{1\/z_j}\rt), \qquad j\in\Bbb Z\setminus\{0\}.
$$
Moreover, we have the identity
$$
\sqrt{\l^2-d^2}={d\/2}\rt(z-{1\/z}\rt)
$$
for all $\l\in \L$ and $z\in \dD$.
If  $\l$ is the eigenvalue of $H$, then we have the identity
\[
\lb{i2}
\begin{aligned}
{d\/2}\rt({1\/z}-z\rt)=-|\l^2-d^2|^{1\/2}\qq if \qq \l<-d,\\
{d\/2}\rt({1\/z}-z\rt)=|\l^2-d^2|^{1\/2}\qq if \qq \l>d.
\end{aligned}
\]

\bigskip
\noindent
The next result follows immediately from Theorem \ref{T4}.
\begin{theorem}
\lb{TrV} {\bf (The trace formulas.)}  Let a real potential $V$ satisfy
\er{Vc}. Then there is infinite number of trace formulae
\[
\lb{rV1} 0\le {\s(\T)\/2\pi}-B_0={1\/2\pi}\int_\T \log |D(e^{it})|dt\le C(d,p)\|V\|_q,
\]
$$
 -\Tr \, V+\sum_{j=1}^N |\l_j^2-d^2|^{1\/2}\sign \l_j
={d\/2\pi}\int_\T e^{-it}d\m(t),
$$
$$
 -\Tr \, V^2+\sum_{j=1}^N  |\l_j||\l_j^2-d^2|^{1\/2}
={d^2\/4\pi}\int_\T e^{-i2t}d\m(t), \qquad \dots.
$$
\end{theorem}


{\bf Proof.} The eigenvalue of $H$, then we have the identity
$$
z={1\/d}\rt(\l \pm\sqrt{\l^2-d^2}  \rt).
$$
\BBox

{\bf Remark.} 1) We consider the case \er{rV1}. If $\l>d$, then we
have $z\in (0,1)$ and then
$$
\begin{aligned}
1-z=1-{1\/d}\rt(\l-\sqrt{\l^2-d^2}  \rt)={1\/d}
\rt(d-\l+\sqrt{\l^2-d^2}\rt)\\
={\sqrt{\l-d}\/d}  \rt(\sqrt{\l+d}-\sqrt{\l-d}\rt)= {2\sqrt{\l-d}\/
(\sqrt{\l+d}+\sqrt{\l-d})}\ge {\sqrt{\l-d}\/ \sqrt{\l+d}}.
\end{aligned}
$$
This yields
$$
\sum_{\l_j>d} {\sqrt{\l_j-d}\/ \sqrt{\l_j+d}}+ \sum_{\l_j<-d}
{\sqrt{-\l_j-d}\/ \sqrt{-\l_j+d}}=\sum_{\l_j} {\sqrt{|\l_j|-d}\/
\sqrt{|\l_j|+d}}\le C_d \|V\|_{2/3}.
$$

\begin{corollary}
\lb{Ter} Let a potential $V$ be real and  satisfy \er{Vc}. Then the
following estimates hold true:
$$
\sum_{j=1}^N  |\l_j||\l_j^2-d^2|^{1\/2}\le \Tr V^2+
{d^2\/4\pi}C(d,p)\|V\|_{2/3}.
$$
\end{corollary}

\bigskip

\section {Appendix, Hardy spaces}
\setcounter{equation}{0}



\subsection{Analytic functions}

We recall the basic facts about the Blaschke product (see pages
53-55 in \cite{G81}) of zeros $\{z_n\}$. The subharmonic function
$v(z)$ on $\O$ has a harmonic majorant if there is a harmonic
function $U(z)$ such that $v(z) \le U(z)$ throughout $\O$.

We need the following well-known results, see e.g. Sect. 2 from \cite{G81}.


\begin{lemma}
\label{TF2}
 Let $\{z_j\}$ be a sequence of points in $\dD\sm \{0\}$ such that $\sum (1-|z_j|)<\iy$
and let  $m\ge 0$ be an integer . Then  the Blaschke product
$$
B(z)=z^m \prod_{z_j\ne 0}{
|z_j|\/z_j}\rt(\frac{z_j-z}{1-\ol z_j z}\rt),
$$
converges in $\dD$. Moreover, the function $B$ is in $\mH_\iy$ and
zeros of $B$ are precisely the points $z_j$, according to the
multiplicity. Moreover,
\[
\lb{BL3}
\bea |B(z)|\le 1 \qqq \forall \ z\in \dD,
\ena
\]
$$
\lim _{r\to 1} |B(re^{i\vt})|=|B(e^{i\vt})|=1 \qqq   \
almost\ everywhere,\qq \vt\in \T,
$$
$$
\lim _{r\to 1}\int_0^{2\pi}\log |B(re^{i\vt})|d\vt=0.
$$
\end{lemma}



\bigskip
\noindent
Let us recall a well-known result concerning analytic functions in the unit disc, e.g., see Koosis page 67 in
\cite{Koo98}.

\begin{theorem}
\lb{TA1} Let $f$ be analytic in the unit disc $\dD$ and let $z_j\ne
0, j=1,2,..., N\le \iy$ be its zeros labeled by
$$
0<|z_1|\le ...\le |z_j|\le |z_{j+1}\le |z_{j+2}|\le ....
$$
Suppose that $f$ satisfies the condition
$$
\sup_{r\in (0,1)}
\int_0^{2\pi}\log |f(re^{i\vt})|d\vt<\iy.
$$
Then
$$
\sum _{j=1}^N (1-|z_j|)<\iy.
$$
The Blaschke product $B(z)$ given by
$$
B(z)=z^m\prod_{j=1}^N {|z_j|\/z_j}{(z_j-z)\/(1-\ol z_j z)},
$$
where $m$ is the multiplicity of $B$ at zero,
 converges absolutely for $\{|z|<1\}$.
Besides,
$f_B(z)=f(z)/B(z)$ is analytic in the unit disc $\dD$ and has no
zeros in $\dD$.

\noindent
Moreover,  if $f(0)\ne 0$ and if $u(z)$ is the least harmonic
majorant of $\log |f(z)|$, then
$$
\sum (1-|z_j|)<u(0) - \log | f (0)|.
$$
\end{theorem}



\bigskip
We now consider  the canonical representation \er{cr} for a
function $f\in \mH_p, p>0$ (see, \cite{Koo98}, p. 76).

\begin{lemma}
\label{Tft}  Let a function $f\in \mH_p, p>0$. Let $B$ be its
Blaschke product. Then there exists a singular measure $\s=\s_f\ge
0$ on $[-\pi,\pi]$ with
\[
\lb{cr}
\begin{aligned}
 f(z)=B(z)e^{ic-K_\s (z)}e^{K_f(z)},\\
K_\s(z)={1\/2\pi}\int_{-\pi}^{\pi}{e^{it}+z\/e^{it}-z}d\s(t),\\
K_f(z)= {1\/2\pi}\int_{-\pi}^{\pi}{e^{it}+z\/e^{it}-z}\log
|f(e^{it})|dt,
 \end{aligned}
\]
for all $|z|<1$, where $c$ is real constant and $\log |f(e^{it})|\in L^1(-\pi,\pi)$.
\end{lemma}

\bigskip
\noindent
We define the functions (after Beurling) in the disc by
$$
\begin{aligned}
& f_{in}(z)=B(z)e^{ic-K_\s (z)} \qq & the \ inner \ factor\ of\ f,\\
&f_{out}(z)=e^{K_f(z)} \qq & the \ outer \ factor\ of\ f,\\
 \end{aligned}
$$
for $|z|<1$. Note that we have  $| f_{in}(z)|\le 1$, since $d\s\ge
0.$

\noindent
Thus $f_B(z)={f(z)\/B(z)} $ has no zeros in the disc $\dD$ and satisfies
$$
\log f_B(z)=ic+{1\/2\pi}\int_{-\pi}^{\pi}{e^{it}+z\/e^{it}-z}d\m(t),
$$
where the measure $\m$ equals
$$
 d\m(t)=\log
|f(e^{it})|dt-d\s(t).
$$

\medskip
\noindent
For a function $f$ continuous on the disc $\ol\dD$ we define the set
of zeros of $f$ lying on the boundary $\pa \dD $ by
$$
\gS_0(f)=\{z\in \dS: f(z)=0\}.
$$
It is well known that the support of the singular measure $\s=\s_f$
satisfies
\[
 \lb{meraze}
\supp \s_f\ss \gS_0(f)=\{z\in \dS: f(z)=0\}
\]
 see for example, Hoffman \cite{Ho62},  p. 70.

\bigskip
\noindent
In the next statement we present trace formulae  for a function $f\in \mH_p, p>0$.

\begin{lemma}
\label{TAt} Let $f\in \mH_p, p>0$ and $f(0)=1$ and let
$B$ be its Blaschke product. Let the functions $\log f$ and $F=\log
f_B$ have the Taylor series in some small disc $\dD_r, r>0$ given by
\[
\lb{asf1}
\begin{aligned}
\log f(z)=-f_1z-f_2z^2-f_3z^3-.....,\\
F=\log f_B(z)=F_0+F_1z+F_2z^2+F_3z^3+.....,\\
\log  B(z)=B_0-B_1z-B_2z^2-..., \qqq as \qqq z\to 0,\\
 F_0=-\log B(0)>0,\qqq F_n=B_n-f_n,\qqq n\ge 1.
\end{aligned}
\]
 Then  the factorization \er{cr} holds true and we have
\[
\lb{ftr0} c=0,\qqq F_0=-\log B(0)={\m(\T)\/2\pi}\ge 0,\qqq 
\]
\[
\lb{ftr1} F_n={1\/\pi}\int_{-\pi}^{\pi}e^{-int}d\m(t),\qqq
n=1,2,....,
\]
where the measure $d\m(t)=\log |f(e^{it})|dt-d\s(t)$.
\end{lemma}

\bigskip
\noindent
{\bf Proof.} Recall that the identity \er{cr} gives
$f(z)=B(z)e^{ic-K_\s (z)}e^{K_f(z)}$, then at $z=0$ we obtain
$$
1=f(0)=B(0)e^{ic-K_\s (0)}e^{K_f(0)}.
$$
Since $B(0), K_\s (0), K_f(0)$ and $c$ are real  we obtain $c=0$.
Moreover, the inequality \eqref{BL3} implies $F_0\ge0$.

\noindent
In order to show \er{ftr1} we need the asymptotics
of the Schwatz integral
\[
\lb{Si} f(z)=B(z)f_B(z),\qqq
 F(z)=\log f_B(z)={1\/2\pi}\int_{-\pi}^{\pi}{e^{it}+z\/e^{it}-z}d\m(t),
\]
 as $z\to 0$. The following identity holds true
\[
\lb{ts1}
{e^{it}+z\/e^{it}-z}=1+{2ze^{-it}\/1-ze^{-it}}=1+2\sum_{n\ge 1}
\big({ze^{-it}}\big)^n=
1+2\big({ze^{-it}}\big)+2\big({ze^{-it}}\big)^2+.....
\]
Thus \er{Si}, \er{ts1} yield the Taylor series at $z=0$:
\[
\lb{asm}
{1\/2\pi}\int_{-\pi}^{\pi}{e^{it}+z\/e^{it}-z}d\m(t)={\m(\T)\/2\pi}+\m_1z+
\m_2z^2+\m_3z^3+\m_4z^4+...\qqq as \qqq z\to 0,
\]
where
$$
\m(\T)=\int_0^{2\pi}d\m(t),\qqq \m_n={1\/\pi}\int_0^{2\pi}e^{-in\vt}
d\m(t),\qqq n\in \Z.
$$
 Thus comparing \er{asf1} and \er{asm} we obtain
 $$
-\log B(0)={\m(\T)\/2\pi}\ge 0,\ F_n=\m_n \qqq \forall \  n\ge 1.
$$
 \BBox


\section {Appendix, estimates involving Bessel's functions}
\setcounter{equation}{0}

\medskip
\noindent
In order complete the proof of Theorem \ref{T2x} we need some uniform estimates for the Bessel functions $J_m, m\in \Z$ with respect to $m$ for which we their integral representation
\begin{equation}
\lb{Be1} J_m(t)={1\/2\pi}\int_0^{2\pi} e^{-imk-i{t}\sin
k}\,dk={i^m\/2\pi}\int_0^{2\pi} e^{-imk+i{t}\cos k}\,dk.
\end{equation}
Note that for all $(t,m)\in \R\ts \Z$:
\begin{equation}
\lb{Be2}
J_{-m}(t)=J_m(t) \qq
\textup{and} \qq  J_{m}(-t)=(-1)^mJ_n(t).
\end{equation}
Our estimates are based on the following three asymptotic formulae, see \cite{Smirn}, Ch IV, $\S$ 2.
Let
$$
\xi = \frac{m}{t}.
$$
Then for a fixed $\varepsilon$, $0<\varepsilon<1$ we have

{\bf 1.} if $\xi>1+\varepsilon$ then
\begin{equation}\label{bess1}
J_m(t) = \frac12 \, \sqrt{\frac{2}{\pi t}} \, \frac{1}{(\xi^2-1)^{1/4}} \, e^{-t\left(\xi \ln\left(\xi + \sqrt{\xi^2 -1}\right) - \sqrt{\xi^2-1}\right)} \, \left(1 + O\left(\frac{1}{m}\right)\right).
\end{equation}
Therefore if $\xi>1+\varepsilon$ then this formula implies the uniform with respect to $m$ exponential decay of the Bessel function in $t$.

\medskip
{\bf 2.} if $\xi<1-\varepsilon$ then
\begin{multline}\label{bess2}
J_m(t) = \frac12\,  \sqrt{\frac{2}{\pi t}} \, \frac{1}{(\xi^2-1)^{1/4}} \left( e^{-i\pi/4i + it\left(-\xi \arccos \xi + \sqrt{1-\xi^2}\right)}
+ e^{i\pi/4i + it\left(\xi \arccos \xi - \sqrt{1-\xi^2}\right)}\right) \\
\left( 1 + O\left(\frac{1}{t}\right)\right).
\end{multline}
In this case the Bessel function oscillates as $t\to\infty$ and obviously the latter formula implies the uniform with respect to $m$ estimate
$$
|J_m(t)|Ê\le C \, t^{-1/2}, \qquad C= C(\varepsilon).
$$

\medskip
{\bf 3.} We now consider the third case $ 1-\varepsilon\le \xi \le 1+ \varepsilon$ which is more difficult.


\begin{lemma}\label{case3}
If $1-\varepsilon\le \xi \le 1+ \varepsilon$, $\varepsilon>0$,  then there is a constant $C=C(\varepsilon)$ such that
\begin{equation}\label{3case}
J_m(t) \le C \, t^{-\frac14} \, \left(|t|^{\frac13}Ê+ |m-t|\right)^{-\frac14}, \qquad \forall \, m, \, t, \quad |m-t|<\varepsilon\, t.
\end{equation}
\end{lemma}

\medskip
\noindent
{\bf Proof.} If $1-\varepsilon\le \xi \le 1+ \varepsilon$, then (see \cite{Smirn}, Ch IV, $\S$ 2)
\begin{equation}\label{bess3}
J_m(t) = \frac{v\left( t^{2/3} \tau(\xi)\right) }{t^{1/3}} \left(c_0(\xi) + O\left(\frac{1}{t}\right)\right) \\+
 \frac{v'\left( t^{2/3} \tau(\xi)\right)}{t^{4/3}} \left(d_0(\xi) + O\left(\frac{1}{t}\right)\right),
\end{equation}
where $v$ is the Airy function and
$$
\tau^{3/2}(\xi) = \xi \, \ln\left(\xi + \sqrt{\xi^2 -1} \right) - \sqrt{\xi^2 -1}
$$
and therefore
\begin{equation}\label{tau}
\tau(\xi) = 2^{1/3}(\xi-1) + O\left((\xi-1)^2\right), \quad {\rm as} \quad \xi \to 1.
\end{equation}
Besides the functions $c_0(\xi)$ and $d_0(\xi)$ are bounded with respect to $\xi$ and, for example,
$$
c_0(\xi) = \sqrt{\frac{2}{\pi}} \, \left(\frac{\tau(\xi)}{\xi^2-1}\right)^{1/4},
$$
(see \cite{Olver} formulae (10.06), (10.07))

In what follows all the constants depend on $\varepsilon$, $0<\varepsilon<1$, but not on $m$ and $t$.

\noindent
Due \eqref{tau} there are constants $c$ and $C$ such that
\begin{equation}\label{y1}
c(\varepsilon) (1+ |y|)^{\frac14} \le  \left(t^{\frac13}  + |m-t|\right)^{\frac14}\, t^{-\frac{1}{12}}\le
C(\varepsilon) (1+ |y|)^{\frac14},
\end{equation}
where $y = t^{\frac23} \, \tau(\xi)$. Moreover, since $|\xi-1| = |\frac{m}{t} -1 |Ê\le \varepsilon$ we also have
\begin{equation}\label{y2}
 \left(t^{\frac13}  + |m-t|\right)^{\frac14}\, t^{-\frac{1}{12}}\le C(\varepsilon) \,  t^{\frac16}.
\end{equation}
Applying estimates for the Airy functions in \eqref{bess3}
$$
|v(y)| \le C \, (1 + |y|)^{-1/4}, \qquad |v'(y)| \le C \, (1 + |y|)^{1/4}
$$
and using \eqref{y1}, \eqref{y2}  we find that if $|t|\ge1$
\begin{multline*}
|J_m(t)| \le C\left( \frac{1}{(1+|y|)^{\frac14} \, |t|^{\frac13}} + \frac{(1+|y|)^{\frac14}}{|t|^{\frac43}} \right)\\
\le
C\, \left( \frac{t^{\frac{1}{12}}}{t^{\frac13}\, \left(t^{\frac13}Ê+ |m-t|\right)^{\frac14}}Ê
+ \frac{\left(t^{\frac13} + |m-t|\right)^{\frac14}}{t^{\frac43}  \, t^{\frac{1}{12}}} \right)\\
\le C \, t^{-\frac14} \, \left(|t|^{\frac13}Ê+ |m-t|\right)^{-\frac14}.
\end{multline*}
The proof is complete. \BBox


\bigskip
\noindent
 Let us now consider the operator $e^{it \D}, t\in \R$. It is unitary  on $L^2(\Bbb S^d)$ and its kernel
is given by
\[
\lb{kr1} (e^{it \D})(n-n')={1\/(2\pi)^{d}}\int_{\Bbb S^d}
e^{-i(n-n',k)+ith(k)}dk,\qqq n,n'\in \Z^d.
\]
where $h(k)=\sum_{j=1}^d\cos k_j, k=(k_j)_{j=1}^d\in \Bbb S^d$.


\begin{lemma}\lb{Texp1}
Let $n=(n_j)_{j=1}^d\in \Z^d, d\ge 1$. Then
\[
\lb{ehtd} (e^{it\D})(n)=i^{-|n|}\prod_{j=1}^d J_{n_j}(t),
\qq
 (n,t)\in \Z^d\times \R,
\]
where $|n|=|n_1|+....+|n_d|$. Moreover, the following estimates are satisfied:
\[
\lb{R12} |(e^{it\D})(n)|\le   C_1 |t|^{-{d\/3}}, \qquad t\ge1,
\]
\[
\lb{iJm} \int_1^\iy |J_{m}(t)|^d\, dt<C_2,\qqq if \qq m\in \Z,\ d\ge 3,
\]
for all  $ (t,n)\in \R\ts \Z^d$ and some constants
$C_1=C_1(d)$ and $C_2=C_2(d)$.
\end{lemma}

\medskip
\noindent
{\bf Proof.} Let $d=1$ and $h(k)=\cos k, k\in \T$. Then  using
\er{kr1}   and  \er{Be1} we obtain
$$
(e^{it\D})(n)={1\/(2\pi)}\int_{\T} e^{-ink+it\cos
k}dk=i^{-n}J_n(t),\qq \forall \ (n,t)\in \Z\ts \R.
$$
which yields \er{ehtd} for $d=1$. Due to the separation of variables we also obtain \er{ehtd} for
any $d\ge 1$.

\noindent
In view of \er{Be2} it is enough to consider the case $n_j\ge 0$.
In order to obtain \eqref{R12} it is enough ti apply the inequalities
\eqref{bess1},  \eqref{bess2} and also  \eqref{3case} if in this inequality we ignore the term $|m-t|$.

\noindent
If $d>3$, then  \eqref{R12} implies \eqref{iJm}. Let now $d=3$. From \eqref{3case} we obtain
\begin{multline*}
 \int_1^\iy |J_{m}(t)|^d\, dt \le C \left(  \int_1^\iy  |t|^{-3/2}\, dt +  \int_1^\iy |t|^{-{3\/4}}\rt(|t|^{1\/3}+|m-t|\rt)^{-{3\/4}} \, dt \right)\\
 \le  C/2 +  C \, \int_1^\iy |t|^{-{3\/4}}\, (1 + |m-t|)^{-{3\/4}}\, dt\\
  \le C/2 + C \, \left( \int_1^\iy |t|^{-3/2}\, dt\right)^{1/2}
  \left( \int_1^\iy (1 + |m-t|)^{-3/2}\, dt\right)^{1/2}<\infty.
\end{multline*}
\BBox


\begin{theorem}
\lb{TApp} i) Let $d\ge 3$.
  Then for each $n\in \Z^d$ the following
estimate holds true:
\[
\lb{Int} \int_1^\iy \bigl|(e^{\pm it\D})(n)\bigr|\, dt\le \b, \qq
\]
where
\[
\lb{Intb} \b=\sup_{m\in \Z}\int_1^\iy |J_{m}(t)|^ddt<\iy.
\]
ii) Let a function $q\in \ell^2(\Z^d)$ and let $X(\l) =qR_0(\l)q,\
\l\in \L$. Then the operator-valued function $X:\L\to \cB_2$ ia
analytic and satisfies
\begin{equation}
\lb{X2} \sup_{\l\in \L} \|X(\l)\|_{\cB_2}\le (1+\b)\|q\|_{2}^2 ,
\end{equation}
\end{theorem}


\medskip
\noindent
{\bf Proof.}

\noindent
i)  Note that \er{iJm} gives \er{Intb}. Due to \er{Be2}
it is sufficient to show \er{Int} for $n\in (\Bbb Z_+)^d$.
Using \er{ehtd} and \er{Intb}, we obtain
$$
\int_1^\iy |(e^{it\D})(n)|dt=\int_1^\iy  \prod_1^d |J_{n_j}(t)|dt\le
\prod_1^d \rt (\int_1^\iy |J_{n_j}(t)|^d\rt)^{1/d}\le  \b,
$$
which yields \er{Int}.

\medskip
\noindent
 ii) Consider the case $\C_-$, the proof for $\C_+$ is
similar. We have the standard representation of the free resolvent
$R_0(\l)$ in the lower  half-plane $\C_-$ given by
$$
\begin{aligned}
R_0(\l)=-i\int_0^\iy e^{it(\D-\l)}dt=R_{01}(\l)+R_{02}(\l),
\\
R_{01}(\l)=-i\int_0^1 e^{it(\D-\l)}dt,\qqq R_{02}(\l)=-i\int_1^\iy
e^{it(\D-\l)}dt,
\end{aligned}
$$
for all $\l\in \C_-$.
 Here the operator valued-function
$R_{01}(\l)$ has analytic extension from  $\C_-$ into the whole
complex plane $\C$ and  satisfies
$$
\|R_{01}(\l)\|\le 1,\qqq \qqq \|qR_{01}(\l)q\|_{\cB_2}\le
\|q\|_{2}^2\qq \forall \ \l\in \C_-,
$$
Let $R_{02}(n'-n,\l)$ be the kernel of the operator $R_{02}(\l)$. We have the identity
$$
R_{02}(m,\l)=-i\int_1^\iy (e^{it(\D-\l)})(m)\,dt,\qqq m=n'-n.
$$
Then the estimate \er{Int} gives
$$
|R_{02}(m,\l)|\le \int_1^\iy|(e^{it\D})(m)|dt \le \b,
$$
which yields
$$
\|qR_{02}(\l)q\|_{\cB_2}^2=\sum_{n,n'\in\Z^d}|q(n)|^2
|R_{02}(n-n'\l)|^2|q(n)|^2\le \sum_{n,n'\in\Z^d}|q(n)|^2
\b^2|q(n)|^2= \b^2\|q\|_{2}^4,
$$
and summing results for $R_{01}$ and $R_{02}$ we obtain \er{X2}.
\BBox


\medskip

\noindent
\textbf{Acknowledgments.}  Various
parts of this paper were written during Evgeny Korotyaev's stay  in
KTH and Mittag-Leffler Institute, Stockholm. He is grateful to the
institutes for the hospitality. He is also grateful to Alexei
Alexandrov (St. Petersburg) and Konstantin Dyakonov (Barcelona),
Nikolay Shirokov (St. Petersburg) for stimulating discussions and
useful comments about Hardy spaces. Our study was supported by the
RSF grant No 15-11-30007.


\end{document}